\documentclass[12pt]{article}
\usepackage{amsfonts}

\usepackage{latexsym}
\usepackage{amssymb}
\usepackage{epsfig}
\usepackage{amsmath}
\usepackage{amsthm}
\usepackage[mathscr]{eucal}
\usepackage{subfigure}
\usepackage{graphicx}

\newtheorem{Lemma}{Lemma}

\newtheorem{Theorem}[Lemma]{Theorem}

\renewcommand{\qed}{\hfill{\ \ \rule{2mm}{2mm}} \vspace{0.2in}}

\newcommand{\ind}{1\hspace{-2.3mm}{1}}

\begin{document}

\title{Phase transitions for Erd\(\ddot{o}\)s-R\(\acute{e}\)nyi graphs}
\author{ \textbf{Ghurumuruhan Ganesan}
\thanks{E-Mail: \texttt{gganesan82@gmail.com} } \\
\ \\
EPFL, Lausanne }
\date{}
\maketitle

\begin{abstract}
Consider the complete graph on \(n\) vertices where each edge is independently open with probability \(p,\) or closed otherwise. Phase transitions for such graphs for \(p = \frac{C}{n}\) have previously been studied using techniques like branching processes and random walks. In this paper, we use an alternate component counting argument for establishing phase transition and obtaining estimates on the sum size of the non-giant components. As a corollary, we also obtain estimates on the size of the giant component for \(C\) large: If \(C\) is sufficiently large, there is a positive constant \(M_0  = M_0(C)\) so that with probability at least \(1-e^{-C/100},\) there is a giant component containing at least \(n - ne^{-C/8}\) vertices and every other component contains less than \(M_0 \log{n}\) vertices.

\vspace{0.1in} \noindent \textbf{Key words:} Erd\(\ddot{o}\)s-R\(\acute{e}\)nyi random graphs, phase transition, tree counting argument.

\vspace{0.1in} \noindent \textbf{AMS 2000 Subject Classification:} Primary:
60J10, 60K35; Secondary: 60C05, 62E10, 90B15, 91D30.
\end{abstract}

\bigskip

\section{Introduction} \label{intro}
Consider \(n\) labelled points \(\{1,2,...,n\}\) and for each \(i,j\) join vertices \(i\) and \(j\) by an edge with probability \(p,\) independently of other edges. The resulting random graph \(G\) is also called the Erd\(\ddot{o}\)s-R\(\acute{e}\)nyi graph (Janson et al (2000)). The phase transitions for such graphs for \(p = \frac{C}{n}\) where \(C\) is a positive constant, have been studied (Janson et al (2000), Durrett (2007), Hofstad (2014)) using various techniques like random walks and branching processes. In this paper, we use a combinatorial component counting argument as an alternate method for establishing phase transitions in \(G.\) Constants mentioned throughout are independent of \(n.\)

We have the following result.
\begin{Theorem}\label{th1} Let \(p = \frac{C}{n}\) where \(C\) is a positive constant.\\
(i) If \(C > e^3,\) there are positive constants \(M_1\) and \(\epsilon_1\) so that \[\mathbb{P}(E(M_1)) \geq \epsilon_1\] for all \(n\) sufficiently large, where \(E(M)\) denotes the event that there is a component of \(G\) containing at least \(\frac{n}{2}\) vertices and every other component contains less than \(M\log{n}\) vertices.\\\\
(ii) If \(C < e^{-1},\) there is a positive constant \(M_2\) so that \[\mathbb{P}(F(M_2)) \geq 1 - \frac{1}{n^2},\] for all \(n\) sufficiently large, where \(F(M)\) denotes the event that every component of \(G\) contains less than \(M\log{n}\) vertices.
\end{Theorem}

For large \(C\) we also have estimates on the sum size of the non-giant components.
\begin{Theorem}\label{th2} There is a positive constant \(C_0\) so that the following holds: If \(C > C_0,\) there is a positive constant \(M = M(C)\) so that \[\mathbb{P}(E(M)) \geq 1-e^{-\frac{C}{100}}\] for all \(n\) sufficiently large, where \(E(M)\) denotes the event that there is an open component of \(G\) containing at least \(n - ne^{-C/8}\) vertices and every other component contains less than \(M\log{n}\) vertices.
\end{Theorem}

The paper is organized as follows: In Section~\ref{pf1}, we prove Theorems~\ref{th1} and Theorem~\ref{th2}.

\section{Proof of Theorem~\ref{th1} and Theorem~\ref{th2}} \label{pf1}
\emph{Proof of (i)}: The first step is to prove that with positive probability there exists a component with size at least \(\frac{n}{2}.\) The second step is to prove that with high probability other components have size less than \(M_0\log{n}\) for some large constant \(M_0 >0.\)

\underline{Step 1}: Let \(C_1\) denote the component containing the vertex \(1.\) If \(C_1\) contains \(r\) vertices, there is a tree \({\cal T} \subset C_1\) containing \(r-1\) edges, each of which is open. Moreover,  every edge with one endvertex in \({\cal T}\) and one endvertex outside, is closed. The number of such edges is \(r(n-r).\) Thus for a fixed tree \({\cal T}_1\) containing \(r\) vertices, we have that \begin{equation}
\mathbb{P}\left(\{\#C_1 = r\}\cap \{C_1 \text{ contains } {\cal T}_1\}\right) \leq p^{r-1}(1-p)^{r(n-r)}. \label{eq_c_bg}
\end{equation}
The number of possible choices for \({\cal T}_1\) is \(r^{r-2}\) and the number of ways of choosing \(r-1\) vertices (apart from the vertex \(1\)) is \({n \choose r-1}.\) Therefore, we have that
\begin{eqnarray}
\mathbb{P}(\#C_1 = r) &\leq& {n \choose r-1} p^{r-1}(1-p)^{r(n-r)}r^{r-2} \label{eq_c_r}
\end{eqnarray}
Alon and Spencer (2003) use a slightly different tree estimate for comparison with a Poisson branching process.

For \(r = 1,\) we have
\begin{equation}
\mathbb{P}(\#C_1 = 1) = (1-p)^{n-1} \leq e^{-C + Cn^{-1}} \label{eq_c_1}
\end{equation}
For \(r \geq 2,\) we have
\begin{eqnarray}
\mathbb{P}(\#C_1 = r) &\leq& \left(\frac{ne}{r-1}\right)^{r-1} p^{r-1} \exp\left(-\frac{Cr(n-r)}{n}\right) \nonumber\\
&=& e^{-C}\left(Ce e^{-C}\right)^{r-1} \left(1-\frac{1}{r}\right)^{-(r-1)} \frac{1}{r}\exp\left(\frac{Cr^2}{n}\right) \nonumber\\
&\leq& \frac{1}{Cr} \left(Ce e^{-C}\right)^{r}\exp\left(\frac{Cr^2}{n}\right) \nonumber\\
&=& \frac{1}{Cr} \exp\left(-r\left\{C- 1 - \log{C} -\frac{Cr}{n}\right\}\right) \label{eq_c_2}
\end{eqnarray}
where we use the estimates \({n \choose k} \leq \left(\frac{ne}{k}\right)^{k}\) and \(1-x \leq e^{-x}\) for the first inequality and the fact that \(\left(1-\frac{1}{r}\right)^{r-1} \geq e\) for all \(r \geq 2\) for the next inequality.

For \(\theta  = \frac{1}{2} + 10^{-3}\) and \(r \leq \theta n,\) we have
\begin{equation}
C- 1 - \log{C} -\frac{Cr}{n} \geq C (1-\theta) -1 - \log{C} =: \delta > \log{2} \label{eq_del}
\end{equation}
for all \(C>e^{3}.\) Indeed, it suffices to check whether \[g(C) := C(1-\theta) - 1 - \log{C} -\log{2} > 0\] for all \(C > e^{3} .\) Since \(e^{3} > \frac{1}{1-\theta}\) and the function \(g(x)\) is increasing for all \(x > \frac{1}{1-\theta},\) we check if \(g(e^{3}) > 0.\) This is true since \(e^3(1-\theta) > 16.\frac{1}{2} - \frac{27}{1000} >7\) and \(1 + \log\left(e^{3}\right) + \log{2} < 1+3+1 = 5.\)

Therefore we have for \(2 \leq r \leq \theta n\) that \[\mathbb{P}(\#C_1 = r) \leq \frac{1}{Cr}e^{-\delta r}.\] Letting \(X_1 = \#C_1 \ind(1 \leq \#C_1 \leq \theta n),\) we have that
\begin{eqnarray}
\mathbb{E}X_1 &\leq& e^{-C + Cn^{-1}} + \frac{1}{C}\sum_{r \geq 2}e^{-\delta r} \nonumber\\
&=& e^{-C + Cn^{-1}} + \frac{e^{-\delta}}{C(e^{\delta}-1)} \nonumber\\
&\leq& e^{-C(1-\theta)} + \frac{e^{-\delta}}{C(e^{\delta}-1)} \nonumber\\
&=& \frac{1}{C(e^{\delta}-1)} \left(e^{-C(1-\theta)+\log{C}}(e^{\delta}-1) + e^{-\delta}\right) \nonumber\\
&=& \frac{1}{C(e^{\delta}-1)} \left(e^{-1} e^{-\delta}(e^{\delta}-1) + e^{-\delta}\right) \nonumber\\
&\leq& \frac{1}{C(e^{\delta}-1)} \left(e^{-\delta}(e^{\delta}-1) + e^{-\delta}\right) \nonumber\\
&=& \frac{1}{C(e^{\delta}-1)}  =: \gamma \label{eq_gm}
\end{eqnarray}
for all \(n\) sufficiently large. From (\ref{eq_del}), we have that \(e^{\delta} > 2\) for all \(C > e^{3}\) and therefore the final term \(\gamma\)  in (\ref{eq_gm}) is strictly less than one.

Choose \(\alpha  > 0\) small such that \(\gamma(1+\alpha) < 1.\) Using Markov's inequality, we have
\begin{equation}
\mathbb{P}\left(\sum_{i=1}^{n} X_i > \gamma(1+\alpha)n \right) \leq \frac{\sum_{i=1}^{n} \mathbb{E}X_i}{\gamma(1+\alpha)n} \leq \frac{1}{1+\alpha}. \label{eq_sum}
\end{equation}
The sum \(\sum_{i=1}^{n} X_i\) is an upper bound on the sum of sizes of components whose size lies between \(1\) and \(\theta n.\) By a size of the component, we mean the number of vertices in the component. Since \(\gamma(1+\alpha)n  < n\) and the total number of vertices is \(n,\) we therefore we obtain that
\begin{equation}
\mathbb{P}(A_{\theta}) \geq \frac{\alpha}{1+\alpha}, \label{eq_a}
\end{equation}
where \(A_{\theta} := \bigcup_{1 \leq i \leq n}\{\#C_i \geq \theta n\}\) denotes the event that there exists a component with at least \(\theta n\) vertices. This proves the first step in the proof.

\underline{Step 2}: To prove that two large components cannot exist, we first obtain from (\ref{eq_c_2}) and (\ref{eq_del}) that
\[\mathbb{P}(\#C_1 = r) \leq e^{-\delta r}\] for \(r \leq \theta n.\) In particular, for a fixed integer \(M\) and \(M\log{n} \leq r \leq \theta n,\) we have  \[\mathbb{P}(\#C_1 = r) \leq \frac{1}{n^{10}}\] if \(M\) is large. Fixing such an \(M,\) we then have that \[\mathbb{P}\left( M\log{n} \leq \#C_1 \leq \theta n\right) \leq \frac{1}{n^8}\] and therefore that
\begin{equation}
\mathbb{P}(B_{\theta}) \leq \frac{1}{n^9}, \label{eq_b}
\end{equation}
where \[B_{\theta} := \bigcup_{1 \leq i \leq n}\left\{M\log{n} \leq \#C_i \leq \theta n\right\}\] denotes the event that there exists a component with size in the range \([M\log{n}, \theta n].\)

From (\ref{eq_a}) and (\ref{eq_b}) , we obtain that \[\mathbb{P}(A_{\theta} \cap B_{\theta}^c) \geq \frac{\alpha}{1+\alpha} - \frac{1}{n^{8}} \geq \frac{\alpha}{2(1+\alpha)}\] for all \(n\) large. If \(A_{\theta} \cap B_{\theta}^c \) occurs, there is at least one component containing at least \(\theta n \geq \frac{n}{2} + 10\) vertices, for all \(n\) large. This means that every other component has at most \(\frac{n}{2}-10\) vertices. But since \(B_{\theta}^c\) also occurs, every other component has size at most \(M\log{n}.\)

\emph{Proof of (ii)}: From (\ref{eq_c_r}) we have for \(r \geq 2\) that
\begin{equation}
\mathbb{P}(\#C_1 = r) \leq {n \choose r-1} p^{r-1}r^{r-2} \leq \left(\frac{enp}{r-1}\right)^{r-1} r^{r-2}
\end{equation}
where we have used \({n \choose k} \leq \left(\frac{ne}{k}\right)^{k}.\) Substituting \(p = \frac{C}{n},\)   the last term above is
 \begin{equation}
\frac{1}{r} (eC)^{r-1}\left(1-\frac{1}{r}\right)^{-r+1} \leq \frac{1}{Cr}(eC)^{r} = \frac{1}{Cr} e^{-\delta_1 r}, \nonumber
\end{equation}
where \(\delta_1 := \log\left(\frac{1}{eC}\right) > 0,\) since \(eC < 1.\) In the first inequality, we have used \(\left(1-\frac{1}{r}\right)^{-r+1} \leq e\) for all \(r \geq 2.\) For a fixed positive integer \(M\) and \(r \geq M\log{n},\) we then have \[\mathbb{P}(\#C_1 =r) \leq e^{-M\log{n}} \leq \frac{1}{n^{10}}\] if \(M\) is large. Fixing such an \(M,\) we have that \[\mathbb{P}\left(\#C_1 \geq M\log{n}\right) \leq \frac{1}{n^{9}}\] and therefore that \[\mathbb{P}\left(H_{\theta}\right) \leq \frac{1}{n^{8}}\] where \[H_{\theta} := \bigcup_{1 \leq i \leq n} \left\{\#C_i \geq M\log{n}\right\}\] denotes the event that there exists a component of size at least \(M\log{n}.\)
\(\qed\)

\emph{Proof of Theorem~\ref{th2}}: As before, let \(\theta = \frac{1}{2} + 10^{-3}\) and for \(1 \leq i \leq n,\) let \(X_i = \#{\cal C}_i \ind(1 \leq \#{\cal C}_i \leq \theta n).\) From (\ref{eq_gm}), we have that \[\mathbb{E}X_i \leq \frac{1}{C(e^{\delta}-1)}.\] For all \(C\) large, we have from (\ref{eq_del}) that \(\delta \geq \frac{C}{4}.\) Therefore, \[\mathbb{E}X_i \leq \frac{2}{C} e^{-C/4} \leq e^{-C/4}\] for all \(C\) large. Setting \(\alpha = e^{C/8} -1\) in (\ref{eq_sum}), we obtain that \[\mathbb{P}\left(\sum_{i} X_i > ne^{-C/8}\right) \leq e^{-C/8} \] for all \(C\) large. From (\ref{eq_a}), we therefore get that
\begin{equation}
\mathbb{P}(A_{\theta}) \geq 1  - e^{-C/8}, \label{eq_a}
\end{equation}
where \(A_{\theta}\) is as in (\ref{eq_a}). From (\ref{eq_b}), we then get that \[\mathbb{P}(A_{\theta} \cap B_{\theta}^c) \geq 1- e^{-C/8} - \frac{1}{n^{8}} \geq 1-2e^{-C/8}\] for all \(n\) large. 

As in the proof of Theorem~\ref{th1}, we have that since \(A_{\theta} \cap B_{\theta}^c\) occurs, there exists at least one component of size at least \(\theta n \geq \frac{n}{2} +10,\) for all \(n\) large. Therefore every other component has at most \(\frac{n}{2} - 10\) vertices. Since \(B_{\theta}^c\) occurs, every other component is of size at most \(M\log{n}.\) Therefore there is exactly one (giant) component containing at least \(n/2 + 10\) vertices. Again since \(A_{\theta}\) occurs, the sum of sizes of all the small components is at most \(ne^{-C/8}.\) This means that the giant component contains at least \(n - ne^{-C/8}\) vertices. \(\qed\)

\section*{Acknowledgements}
I thank Professors Rahul Roy and Thomas Mountford for crucial comments and for my fellowship.



\setcounter{equation}{0} \setcounter{Lemma}{0} \renewcommand{\theLemma}{II.%
\arabic{Lemma}} \renewcommand{\theequation}{II.\arabic{equation}} %
\setlength{\parindent}{0pt}




%





\bibliographystyle{plain}

\end{document}